\documentclass[12pt]{article}

\usepackage{amsmath}
\usepackage{amssymb}
\usepackage{amsthm}
\usepackage{latexsym}
\usepackage{color}
\usepackage{graphicx}
\usepackage{color}

\DeclareSymbolFont{calletters}{OMS}{cmsy}{m}{n}
\DeclareSymbolFontAlphabet{\mathcal}{calletters}

\def\be{\begin{eqnarray}}
\def\ee{\end{eqnarray}}

\def\b*{\begin{eqnarray*}}
\def\e*{\end{eqnarray*}}

\newtheorem{Theorem}{Theorem}[section]

\newtheorem{Proposition}[Theorem]{Proposition}

\newtheorem{Assumption}[Theorem]{Assumption}

\newtheorem{Lemma}[Theorem]{Lemma}
\newtheorem{Corollary}[Theorem]{Corollary}
\newtheorem{Remark}[Theorem]{Remark}

\usepackage{color}

\def \D{\mathbb{D}}
\def \E{\mathbb{E}}

\def \R{\mathbb{R}}
\def \S{\mathbb{S}}

\def \T{\mathbb{T}}

\def\Gc{{\cal G}}

\def\Mc{{\cal M}}

\def\Wc{{\cal W}}

\def \om{\omega}

\def \eps{\varepsilon}
\def \0{\mathbf{0}}

\def \Xh{\widehat{X}}

\def\x{\times}

\def\1{{\bf 1}}
\def \proof{{\noindent \bf Proof. }}

\def\psih{\widehat{\psi}}

\def\Xt{\widetilde{X}}
\def\It{\widetilde{I}}

\def\psit{\widetilde{\psi}}

\def\Xb{\overline{X}}

\def\Wch{\widehat{\Wc}}

\def\Ih{\widehat{I}}
\def\xb{\bar x}
\def\mt{\tilde m}

\title{Unbiased simulation of Asian options}

\author{
Bruno Bouchard\footnote{CEREMADE, Universit\'e Paris-Dauphine, PSL, CNRS.  bouchard@ceremade.dauphine.fr. }
\and Xiaolu Tan\footnote{Department of Mathematics, the Chinese University of Hong Kong. \texttt{xiaolu.tan@cuhk.edu.hk}. Research supported by Hong Kong RGC General Research Fund (project 14302622) and the Faculty Direct Grant.}
}

\date{\today}

\begin{document}

\maketitle

\abstract{
	We provide an extension of the unbiased simulation method for Stochastic Differential Equations (SDEs) developed in Henry-Labord\`ere et al. [{\it Ann Appl Probab.} 27:6 (2017) 1-37]
	to {a class of} path-dependent {dynamics, pertaining to Asian options}.
	In {our} setting, both the payoff and the SDE coefficients  depend on the (weighted) average of the process {or, more precisely, on the integral of the solution to the SDE against a continuous function with finite variation. In particular, this applies to the numerical resolution of the class of path-dependent PDEs whose regularity, in the sense of Dupire, is studied in Bouchard and Tan [{\it Ann.  I.H.P., to appear}]. }

	\vspace{5mm}


	\vspace{5mm}

}

\section{Introduction}

	Let $d \ge 1$, $T > 0$ and $W$ be a $d$-dimensional Brownian motion,
	$\mu: [0,T] \x \R^d \x \R^d \longrightarrow \R^d$ and $\sigma: [0,T] \x \R^d \x \R^d \longrightarrow \S^d$ be the drift and diffusion coefficients,
	where $\S^d$ denotes the collection of all $d \x d$ dimensional matrices.
	We consider the process $X$ defined as {the} solution {to} the path-dependent SDE
	\begin{equation} \label{eq:PSDE_intro}
		X_t
		=
		x_0
		+
		\int_0^t \mu \big(s, X_s, I_s \big) ds
		+
		\int_0^t \sigma \big(s, X_s, I_s \big) dW_s,
		~\mbox{with}~
		I_t := \int_0^t X_s dA_s,
	\end{equation}
	where $W$ is a $d$-dimensional Brownian motion and $A: [0, T] \longrightarrow \R$ is a (deterministic) continuous function with finite variation.
	Our main objective is to provide a Monte-Carlo simulation method to estimate the expected value
	\be \label{eq:V0_intro}
		V_0 &:=& \E \big[ g(X_T, I_T ) \big],
	\ee
	for some function $g: \R^d \x \R^d  \longrightarrow \R$.

	\vspace{0.5em}
	
	{The} numerical simulation of {the} SDE \eqref{eq:PSDE_intro} is usually based on {a} discrete time {approximation} scheme (e.g. the Euler scheme),  so that the corresponding Monte Carlo estimator for $V_0$ in \eqref{eq:V0_intro} {suffers from} a discretization error, {see e.g.~the 
	seminal work of Talay and Tubaro \cite{TalayTubaro}, as well as  Kloeden and Platen \cite{KP} and Graham and Talay \cite{GrahamTalay} for an overview of various discrete time approximation approaches.}
	
	\vspace{0.5em}
	
	 For a special class of SDEs, {including}  one-dimensional homogeneous SDEs with constant volatility coefficient, exact simulation methods of the marginal distribution of $X_T$ are available,  see e.g.~Beskos and Roberts \cite{bes1}, Beskos, Papaspiliopoulos and Roberts \cite{bes2}, Jourdain and Sbai \cite{Jourdain}, etc. 

	\vspace{0.5em}

	One can also get rid of the discretization error in the Monte Carlo estimator {by appealing to the so-called}  unbiased simulation methods. %
	One possibility is to apply a random level in the multilevel Monte Carlo method of Giles \cite{Giles}  in the case where $(X,A)$ is an It\^o process,  see e.g. Rhee and Glynn \cite{RheeGlynn}.  
	In the uniformly elliptic case, i.e.~$A\equiv 0$ and $g(X_T,I_T)=g(X_T)$, another one consists in multiplying the terminal payoff by a suitable  random variable ${\cal M}_T$, also called Malliavin weight, so as to compensate exactly the expectation bias induced by the use of an Euler scheme $\Xh$ in place of $X$. The definition of ${\cal M}_T$ is based on the parametrix method for PDEs, which ensures that
	\be\label{BKH}
		V_0=\E \big[ g\big( \Xh_T \big) \Mc_T \big].
	\ee
	Since $\Xh$ and  ${\cal M}_T$ can be simulated exactly, this provides an unbiased MC estimator for $V_0$, see e.g.~Bally and Kohatsu-Higa \cite{Bally}, Andersson and Kohatsu-Higa \cite{AKH}, Chen, Frikha and Li \cite{ChenFrikhaLi}, as well as Henry-Labord\`ere, Tan and Touzi \cite{HLTT} and  Doumbia, Oudjane and Warin \cite{DoumbiaOudjaneWarin}, Agarwal and Gobet \cite{AG}.

	\vspace{0.5em}

	In this paper, we {extend the approach of} \cite{DoumbiaOudjaneWarin,HLTT}  to the path-dependent setting \eqref{eq:PSDE_intro}-\eqref{eq:V0_intro}. {Even when $A$ is absolutely continuous, e.g. $A_t=t$,   $(X_t, I_t)_{t \ge 0}$  is a degenerate diffusion process and the estimator proposed in \cite{HLTT}  cannot be applied. The reason is that it relies on the Markovian representation of $V_0$ as a function of $(X_0,I_0)$. Rewriting $V_0$ as the time $0$ value of a path-dependent functional, and building on the approach developed in \cite{BouchardTan} to study the (Dupire's) regularity of the path-depend PDE associated to  \eqref{eq:V0_intro}, we consider a pretty general situation in which $A$ does not need to be absolutely continuous. 	
	We  are able to find the  weight function ${\cal M}_T$ associated to the Dupire's vertical derivatives of the associated value function, and  bound the second moment of the estimator so that the Monte Carlo estimation error can be controlled in the usual way. This requires a structural condition relating the H\"older regularity of the coefficients with the regularity of $A$ in the spirit of \cite{BouchardTan}, see Assumption \ref{assum:mu_sigma} and condition \eqref{eq:cond_kappa} below.}

	\vspace{0.5em}

	The paper is organized as follows.
	In Section \ref{sec:pSDE}, we explain {how to construct our unbiased representation for $V_0$ in \eqref{eq:V0_intro}, under general, abstract,  integrability conditions.	Then, in Section \ref{sec:integrability}, we provide    sufficient conditions that ensure the (square) integrability of our estimator.}
	{Finally}, in Section \ref{sec:Numerics}, we provide some numerical  examples in the context of the pricing of Asian options.
	Technical lemmas are {left to} Appendix \ref{sec:technical_lemmas}.

\section{Unbiased estimator for path-dependent SDEs}
\label{sec:pSDE}

	Let us consider the path-dependent SDE in \eqref{eq:PSDE_intro}
	\begin{equation} \label{eq:PSDE}
		d X_t
		=
		\mu \big(t, X_t, I_t \big) dt
		+
		\sigma \big(t, X_t, I_t \big) dW_t,
		~\mbox{with}~
		I_t := \int_0^t X_s dA_s,
	\end{equation}
	where $A: [0, T] \longrightarrow \R$ is a given continuous function with finite variation, and $W$ is a Brownian motion.
	Our objective is to estimate 
	\be \label{eq:def_V0}
		V_0 &:=& \E \big[ g(X_T, I_T ) \big],
	\ee
	in which $(\mu,\sigma,g) : [0,T]\x \R^d\x \R^d\mapsto \R^d\times \S^d\times \R$ are measurable.

	We   assume that the volatility function $\sigma$ is non-degenerate, and that the above path-dependent SDE \eqref{eq:PSDE} has a unique weak solution. {See  Assumption \ref{assum:smooth} below. }
	
\begin{Assumption} \label{assum:smooth}
		$\mathrm{(i)}$ The coefficient functions $\mu$ and $\sigma$ are both continuous and have at most linear growth,
		the volatility coefficient $\sigma$ is non-degenerate in the sense that there exists some constant $\eps_0 > 0$ such that
		\b*
			\sigma \sigma^{\top}(t,x, \xb) ~\ge~ \eps_0 I_d,
			~~\mbox{for all}~ (t, x, \xb) \in [0,T] \x \R^d  \x \R^d.
		\e*
		Further, the path-dependent SDE \eqref{eq:PSDE} has a unique weak solution, {for} any initial condition $(t, x, \xb) \in [0,T] \x \R^d \x \R^d$.
		
		\vspace{0.5em}

		\noindent $\mathrm{(ii)}$ The function $g$ is continuous and has at most polynomial growth.
		Moreover, {the map} $u: [0,T] \x \R^d \x \R^d \longrightarrow \R$   defined by
		\begin{equation} \label{eq:def_u}
			u(t, x, \xb) 
			~:=~
			\E \big[ g \big(X_T, I_T \big) \big| X_t =x, I_t = \xb \big]
		\end{equation}
		{admits first and second order derivatives} $D_x u(t,x, \xb)$ and $D^2_{xx} u(t, x, \xb)$ {with respect to its second argument, that} are continuous and have at most polynomial growth on $[0,T] \x \R^d \x \R^d$.
	\end{Assumption}
	
	\begin{Remark}
		The existence of the derivatives $D_x u(t,x, \xb)$ and $D^2_{xx} u(t, x, \xb)$ can be obtained by working along the lines of \cite{BouchardTan}.  However, they will not enjoy a uniform polynomial growth unless $g$ is smooth enough. This condition can however easily be removed by adding an approximation argument in our proofs, see Remark \ref{rem:smoothness} and Corollary \ref{coro:approx_repres} below. Importantly,  no differentiability assumption is made with respect to the second space variable $\xb$ (which would require all coefficients to be smooth).
	\end{Remark}

\subsection{The representation formula}
\label{subsec:general_SDE}

\paragraph{An intuitive overview.}
	Before giving the explicit representation formula for the expected value $V_0$,
	let us briefly outline the intuition behind its construction,
	but in a $1$-dimensional simple context with coefficient functions $(\mu, \sigma, g)$ being only functions of $x \in \R$, so that 
	$$
		V_0 = \E \big[ g(X_T) \big], ~\mbox{with}~ dX_t = \mu(X_t) dt + \sigma(X_t) d W_t.
	$$
	In this simplified setting, $V_0$ can be characterized by a parabolic PDE, i.e.~$V_0 = u(0, x_0)$ with $u$ being (the unique) solution to the parabolic PDE
	$$
		\partial_t u(t,x) + \mu(x) {D}_x u(t,x) + \frac12 \sigma^2(x) {D}^2_{xx} u(t,x) = 0,
		~~
		u(T, x) = g(x).
	$$
	{The core idea consists in looking at $u$ as the} solution to the PDE
	$$
		\partial_t u(t,x) + \mu_0 {D}_x u(t,x) + \frac12 \sigma^2_0 {D}^2_{xx} u(t,x) + f(t,x) = 0,
	$$
	with {the source term $f$ defined by}
	\begin{equation} \label{eq:def_f_illustration}
		f(t,x) := \big( \mu(x) - \mu_0 \big) {D}_x u(t,x) + \frac12 \big( \sigma^2(x) - \sigma^2_0 \big)  {D}^2_{xx} u(t,x),
	\end{equation}
	{for some given constants $\mu_0$ and $\sigma_0$.}
	It {then} follows {from the} Feynmann-Kac's formula that
	\begin{align*} 
		u(t,x) =
		\E\Big[ g \big(x+ &\mu_0 (T-t) + \sigma_0 (W_T - W_t) \big)  \nonumber \\
		&+ \int_t^T \!\!  f \big(s, x+ \mu_0(s-t) + \sigma_0 (W_s - W_t) \big)  ds \Big] .
	\end{align*}
	{If $f$ was known, the above could be simulated without bias, up to the discretization of the integral. However, it depends on the derivatives of $u$, which is precisely the unknown we want to estimate, and we also want to get rid of any biais du to the integral term.}
	{As for the integral term, we simply}  introduce a {$\R_+$-valued} random variable $\tau$, which is independent of the Brownian motion $W$, and {with} distribution function $F$ and density function $\rho$.
	Then
	\begin{align} \label{eq:rep_u_psi1}
		u(t,x) = \E \big[ \psi^{t,x}_1 \big],
	\end{align}
	with
	$$
		\psi^{t,x}_1 :=
		\frac{g \big(x + \mu_0 \tau + \sigma_0 (W_{t+ \tau} - W_t) \big)}{1 - F(T-t)} \1_{\{\tau > T\}}
			+ \frac{f \big(\tau, x+ \mu_0 \tau + \sigma_0 (W_{t+\tau} - W_t) \big)}{\rho(\tau)}.
	$$
	{To get rid of the derivatives in $f$, we use the fact that} $u(t,x)$ can be represented as the expected value of a functional of the increments of the Brownian motion in \eqref{eq:rep_u_psi1}. {Thus,}
	one can {use a standard integration by parts argument to rewrite the derivatives of $u$ in $f$ as}
	\begin{equation} \label{eq:pre_du12}
		\partial_x u(t,x) = \E \big[ \psi^{t,x}_1 \Wc_1 \big],
		~~
		\partial^2_{xx} u(t,x) = \E \big[ \psi^{t,x}_1 \Wc_2 \big],
	\end{equation}
	for some (Malliavin) weight functions $\Wc_1$ and $\Wc_2$ {that do not involve $u$}. {One can then plug \eqref{eq:pre_du12} in \eqref{eq:def_f_illustration} to obtain  another representation in terms of a random variable $\psi^{t,x}_2$, which will also depend on $f$ but at a higher order. Then, one plugs \eqref{eq:pre_du12} in $f$  again in the definition of $\psi^{t,x}_2$... Iterating this procedure, together with good choices of $(\mu_0, \sigma_0)$, leads to the definition of a sequence  $(\psi^{t,x}_n)_{n \ge 1}$ of representation variables, whose limit does not depend on $u$ anymore but is composed of a random number of finitely many terms, which can be simulated exactly.}
	The detailed induction procedure, that we adapt here to a path-dependent PDE context, will be presented in the proofs of Section \ref{subsec:proof_mainthm}.

\paragraph{The representation formula.}

	Let us now provide the representation formula for $V_0$ in the path-dependent setting \eqref{eq:def_V0}.

	\vspace{0.5em}
	
	{Let $\rho: \R \longrightarrow \R_+$ be the density function of some distribution supported on $\R_+$, and $F: \R \longrightarrow [0,1]$ be the corresponding c.d.f.~(cumulative distribution function). We assume that $\rho>0$ on $[0,T]$ and that $F(T)<1$.
	We consider a sequence $(\tau_k)_{k \ge 1}$ of i.i.d.~random variables following the distribution with c.d.f.~$F$, {independent of $W$.}
	We then define the random time grid $(T_k)_{k \ge 1}$ by $T_0 := 0$ and 
	\begin{equation} \label{eq:def_Tk}
		T_{k}  := \big( T_{k-1} + \tau_{k} \big) \wedge T.
		~~k \ge 1,
	\end{equation}
	Set
	$$
		\Delta T_{k+1} := T_{k+1} - T_k,~~k \ge 0,
		~~~
		N_T := \max \big\{ k \ge 0 ~: T_k < T \big\}
	$$
	so that $N_T$ denotes the number of points of the grid $(T_k)_{k \ge 0}$ belonging to $[0,T)$.}
	
	{We then}  define the process $(\Xh, \Ih)$ {as} the Euler scheme of \eqref{eq:PSDE} on the random {time} grid $(T_k)_{k \ge 0}$, that is
	$$
		\Xh_0 := x_0,
		~~~\Ih_0 = 0,
	$$
	and then, for each $k = 0, 1, \cdots, N_T$,
	$$
		\Xh_t 
		~:=~
		\Xh_{T_k}
		+
		\mu\big(T_k, \Xh_{T_k}, \Ih_{T_k} \big) \big(t - T_{k} \big)
		+
		\sigma \big(T_k, \Xh_{T_k}, \Ih_{T_k} \big) \big( W_t - W_{T_{k}} \big), 
		~~t \in (T_k, T_{k+1}],
	$$
	and
	\begin{align}\label{eq:def_Ih1t}
		\Ih_t 
		~:=~
		\Ih_{T_k} + \int_{T_k}^t \Xh_s d A_s,
		~~t \in (T_k, T_{k+1}].
	\end{align}
	Equivalently, given $(\Xh_{T_k}, \Ih_{T_k})$, the couple $(\Xh_{T_{k+1}}, \Ih_{T_{k+1}})$ is defined by
	\begin{equation} \label{eq:def_Xh}
		\Xh_{T_{k+1}}
		:=
		\Xh_{T_k}
		+
		\mu\big(T_k, \Xh_{T_k}, \Ih_{T_k} \big) \big(T_{k+1} - T_{k} \big)
		+
		\sigma \big(T_k, \Xh_{T_k}, \Ih_{T_k} \big) \big( W_{T_{k+1}} - W_{T_{k}} \big),
	\end{equation}
	and
	\begin{align} 
		\Ih_{T_{k+1}}
		&{=}
		\Ih_{T_k}
		~+~
		\Xh_{T_k} \big( A_{T_{k+1}} - A_{T_k} \big) 
		+ \mu(T_k, \Xh_{T_k}, \Ih_{T_k} ) \int_{T_k}^{T_{k+1}} (A_{T_{k+1}} - A_s) ds \nonumber \\
		&~~~~~~~~~~~~~~~~~~~~~~~~~~~~~~~~~~~~~~+ \sigma(T_k, \Xh_{T_k}, \Ih_{T_k} ) \int_{T_k}^{T_{k+1}} (A_{T_{k+1}} - A_s) d W_s \label{eq:def_Ih1} \\
		&=
		\Ih_{T_k}
		~+~
		\Xh_{T_{k+1}} \big( A_{T_{k+1}} - A_{T_k} \big) 
		- \mu(T_k, \Xh_{T_k}, \Ih_{T_k} ) \int_{T_k}^{T_{k+1}} (A_s - A_{T_k}) ds \nonumber \\
		&~~~~~~~~~~~~~~~~~~~~~~~~~~~~~~~~~~~~~~- \sigma(T_k, \Xh_{T_k}, \Ih_{T_k} ) \int_{T_k}^{T_{k+1}} (A_s - A_{T_k}) d W_s \label{eq:def_Ih2}.
	\end{align}
	In particular, given $(T_k, T_{k+1})$, the conditional distribution of the random vector 
	$$
		\Big( W_{T_{k+1}} - W_{T_k}, ~\int_{T_k}^{T_{k+1}} (A_s - A_{T_k}) dW_s \Big)
	$$
	is a Gaussian vector,
	so that the sequence $\big(\Xh_{T_k}, \Ih_{T_k} \big)_{k \ge 0}$ can be simulated explicitly. 

	\vspace{0.5em}

	\vspace{0.5em}
	
		From now on let us define, for each $k \ge 0$, 	
	$$
		\mu_{T_k} := \mu \big(T_k, \Xh_{T_k}, \Ih_{T_k} \big),
		~~~
		\sigma_{T_k} := \sigma \big(T_k, \Xh_{T_k}, \Ih_{T_k} \big),
		~~~
		a_{T_k} := \sigma_{T_k} \sigma_{T_k}^{\top},
	$$
	and
	$$
		{D_x u_{T_k} := D_x u \big(T_k, \Xh_{T_k}, \Ih_{T_k} \big),~~~ D^2_{xx} u_{T_k} :=   D^2_{xx} u \big(T_k, \Xh_{T_k}, \Ih_{T_k} \big), }
	$$
	as well as
	$$
		\overline A_{k, k+1} := \frac{1}{\Delta T_{k+1}} \int_{T_k}^{T_{k+1}} A_r dr,
		~~~~
		m^1_{k,k+1} ~:=~ \frac{1}{\Delta T_{k+1}} \int_{T_k}^{T_{k+1}} (A_{r} - A_{T_k}) d r,
	$$
	$$
		m^2_{k,k+1} ~:=~ \frac{1}{\Delta T_{k+1}} \int_{T_k}^{T_{k+1}} \big(A_{r} - \overline A_{k,k+1} \big)^2 d r,
	$$
	and
	$$
		\mt^2_{k,k+1} ~:=~ \frac{1}{\Delta T_{k+1}} \int_{T_k}^{T_{k+1}} (A_{r} - A_{T_k})^2 d r.
	$$
	We also set, for each $k \ge 1$,  
	\be \label{eq:Wcb1}
		\Wch^1_k
		~:=~
		\big(\mu_{T_k} - \mu_{T_{k-1}} \big)
		\cdot~
		M_{k+1},
  	\ee
	and 
	\begin{align} \label{eq:Wcb2}
		\Wch^2_k
		~:=~
		\frac12 \mathrm{Tr} \Big[ \Big(
			& a_{T_k}  - a_{T_{k-1}}
		\Big) 
		\Big( 
			M_{k+1} M_{k+1}^{\top} 
			-
			\frac{1}{\Delta T_{k+1}} \frac{\mt^2_{k, k+1}}{m^2_{k,k+1}} a_{T_k}^{-1}
		\Big)
		\Big],
	\end{align}
	where
	\begin{align} \label{eq:defMk}
		M_{k+1} 
		~:=~
		\frac{1}{\Delta T_{k+1} m^2_{k,k+1} }
		\big( \sigma^{\top}_{T_{k}} \big)^{-1}
		\int_{T_k}^{T_{k+1}}
		\Big(
			\mt^2_{k, {k+1}} 
			-
			m^1_{k, {k+1}}  
			\big( A_s - A_{T_k} \big) 
		\Big) dW_s.
	\end{align}

	Finally, for each $n \ge 1$, using $D_x u$ and $D^2_{xx} u$ defined in Assumption \ref{assum:smooth} 
	and the fact that $\Delta T_{N_T+1 }=T-T_{N_T}$ by construction, we introduce 
	\begin{align} \label{eq:def_psi_n}
		\psih_n 
		&~:=~
			\1_{\{N_T \le n-1\}} ~\frac{g \big(\Xh_T, \Ih_T \big) - g\big(\Xh_{T_{N_T}}, \Ih_{T_{N_T}} \big) \1_{\{N_T > 0\}} }{1 - F(\Delta T_{N_T+1 })} 
				\prod_{k=1}^{N_T} \Big( 
					~\frac{~\Wch^1_k + \Wch^2_k~}{\rho(\Delta T_k)}
				\Big) \nonumber \\
			&+ 
			\1_{\{N_T \ge n\}} ~ \frac{\big(b_{T_n} - b_{T_{n-1}}\big) \!\cdot\! D_x  u_{T_n} + \frac12 \mathrm{Tr} \big[ \big(a_{T_n} - a_{T_{n-1}} \big) D^2_{xx}  u_{T_n} \big] }{~\rho(\Delta T_n)~}~ 
				\prod_{k=1}^{n-1}  \Big( 
					~\frac{~\Wch^1_k + \Wch^2_k~}{\rho(\Delta T_k)}
				\Big)
	\end{align}
	{with the convention that the product over an empty set equals $1$}, 
	and then define the representation random variable $\psih$ by
	\be \label{eq:def_psih}
		\psih
		&:=&
		\frac{g \big(\Xh_T, \Ih_T \big)
		-
		g \big( \Xh_{T_{N_T}}, \Ih_{T_{N_T}} \big) \1_{\{N_T > 0\}}}
		{1 - F(\Delta T_{N_T+1})  }
		~\prod_{k=1}^{N_T} 
		~\frac{~\Wch^1_k + \Wch^2_k~}{\rho(\Delta T_k)}.
	\ee

	Recalling that $T_{N_T} < T$ and $T_{N_T+1} = T$ so that $\Delta W_{N_T+1} = W_T - W_{T_{N_T}}$,
	we introduce $(\Xt_T, \It_T)$ by reversing the sign of the Brownian increment on $[T_{N_T}, T]$ only, i.e.
	$$
		\Xt_T
		:=
		\Xh_{T_{N_T}}
		+
		\mu\big(T_{N_T}, \Xh_{T_{N_T}}, \Ih_{T_{N_T}} \big) \big(T - T_{N_T} \big)
		-
		\sigma \big(T_{N_T}, \Xh_{T_{N_T}}, \Ih_{T_{N_T}} \big) \big( W_{T} - W_{T_{N_T}} \big),
	$$ 
	and
	\begin{align*} 
		\It_T
		&:=
		\Ih_{T_{N_T}}
		~+~
		\Xt_{T} \big( A_{T} - A_{T_{N_T}} \big) 
		- \mu(T_{N_T}, \Xh_{T_{N_T}}, \Ih_{T_{N_T}} ) \int_{T_{N_T}}^{T} (A_s - A_{T_{N_T}}) ds \nonumber \\
		&~~~~~~~~~~~~~~~~~~~~~~~~~~~~~~~~~~~~~~+ \sigma(T_{N_T}, \Xh_{T_{N_T}}, \Ih_{T_{N_T}} ) \int_{T_{N_T}}^{T} (A_s - A_{T_{N_T}}) d W_s.
	\end{align*}
	We then define
	\begin{align*} 
		\psit_n 
		&:=
			\1_{\{N_T \le n-1\}} ~\frac{g \big(\Xt_T, \It_T \big) - g\big(\Xh_{T_{N_T}}, \Ih_{T_{N_T}} \big) \1_{\{N_T > 0\}} }{1 - F(\Delta T_{N_T+1 })} 
				\bigg( \prod_{k=1}^{N_T-1} 
					\frac{~\Wch^1_k + \Wch^2_k~}{\rho(\Delta T_k)}
				\bigg)
				\frac{\Wch^2_{N_T}- \Wch^1_{N_T} }{\rho(\Delta T_{N_T})} \\
			&+ 
			\1_{\{N_T \ge n\}} ~ \frac{\big(b_{T_n} - b_{T_{n-1}}\big) \!\cdot\! D_x  u_{T_n} + \frac12 \mathrm{Tr} \big[ \big(a_{T_n} - a_{T_{n-1}} \big) D^2_{xx}  u_{T_n} \big] }{~\rho(\Delta T_n)~}~ 
				\prod_{k=1}^{n-1}  \Big( 
					~\frac{~\Wch^1_k + \Wch^2_k~}{\rho(\Delta T_k)}
				\Big),
	\end{align*}
	as well as
	\begin{equation*}
		\psit
		~:=~
		\frac{g \big(\Xt_T, \It_T \big)
		-
		g \big( \Xh_{T_{N_T}}, \Ih_{T_{N_T}} \big) \1_{\{N_T > 0\}}}
		{1 - F(\Delta T_{N_T+1})  }
		~\bigg( \prod_{k=1}^{N_T-1} 
		\frac{ \Wch^1_k + \Wch^2_k~}{\rho(\Delta T_k)} \bigg)
		\frac{\Wch^2_{N_T}- \Wch^1_{N_T} }{\rho(\Delta T_{N_T})}.
	\end{equation*}
	Finally, let
	\begin{equation} \label{eq:def_psi}
		\psi_n := \frac{\psih_n + \psit_n}{2}, ~~n\ge 1,
		~~\mbox{and}~~
		\psi := \frac{\psih + \psit}{2}.
	\end{equation}

	\begin{Theorem} \label{theo:renorm_general_vol}
		Let Assumption \ref{assum:smooth} hold true.
		Suppose in addition that $\psi_n$ (resp. $\psih_n$, $\psit_n$) is integrable for each $n \ge 1$,
		then
		$$
			\E \big[ \psi_n \big] = V_0
			~~\mbox{(resp.}~\E \big[ \psih_n \big] = V_0, ~\E \big[ \psit_n \big] = V_0 \mbox{)},
			~~\mbox{for each}~
			n \ge 1.
		$$
		Further assume that the sequence $(\psi_n)_{n \ge 1}$ (resp.  $(\psih_n)_{n \ge 1}$,  $(\psit_n)_{n \ge 1}$) is uniformly integrable,
		then the random variable $\psi$ (resp. $\psih$, $\psit$) is also integrable and
		\begin{equation} \label{eq:unbiased_estim}
			\E \big[ \psi \big] = V_0
			~~\mbox{(resp.}~
			\E \big[ \psih \big] = V_0, ~\E \big[ \psit \big] = V_0
			\mbox{)}.
		\end{equation}
	\end{Theorem}

	\begin{Remark}
		$\mathrm{(i)}$ Note that the random variables $(\psi_n)_{n\ge 1}$ depend on the value function $u$ which is exactly the value we aim to estimate. This sequence is only used to approximate $\psi$ in  the proof of \eqref{eq:unbiased_estim}. 
		The random variable $\psi$ does not depend on $u$ and can be simulated exactly, so that it provides an unbiased square integrable Monte Carlo estimator for $V_0$ by \eqref{eq:unbiased_estim}, as soon as it is square integrable, see  Section \ref{sec:integrability}.
		
		\vspace{0.5em}

		\noindent $\mathrm{(ii)}$ By the symmetry of the distribution of the Brownian motion, the random variable $\psit$ can be considered as {a form of} an antithetic variable of $\psih$. They have the same distribution.
		The reason for considering $\psi$ rather than $\psih$ is that its integrability is easier to handle,
		see more discussions in Remark \ref{rem:antithetique}.

	\end{Remark}

	Before turning to the proof of \eqref{eq:unbiased_estim},  let us also discuss Assumption \ref{assum:smooth}.

	\begin{Remark}
		The non-degeneracy condition on $\sigma$ in Assumption \ref{assum:smooth} is  crucial, {see in particular \eqref{eq:defMk}.}
	\end{Remark}

	\begin{Remark} \label{rem:smoothness} 	The existence, continuity and growth condition on $D_x u$ and $D^2_{xx} u$ in Assumption \ref{assum:smooth} is {only used to define} the random variables $(\psih_n)_{n\ge 1}$ in \eqref{eq:def_psi_n}.
	
	\vspace{0.5em}
	
		\noindent $\mathrm{(i)}$ 
		When $A_t = t$ for all $t \in [0,T]$, and $(\mu, \sigma)$ are bounded and continuous, and satisfy a technical H\"older condition, it is proved in Francesco and Pascucci \cite{FrancPacucci} that $D_x u$ and $D^2_{xx} u$ exist and are continuous. The extension to  a class of general continuous functions $A$ with finite variation was obtained by Bouchard and Tan \cite{BouchardTan}.  {However, a uniform polynomial growth can not be expected unless the coefficients (including the terminal payoff function) are smooth.} 
		
%

		\vspace{0.5em}
		
		\noindent $\mathrm{(ii)}$ The  formulation of  $\psih$   does not depend on $(D_x u, D^2_{xx} u)$. 
		{In particular}, 
		if there exists a sequence of coefficients $(\mu^m, \sigma^m, g^m)_{m \ge 1}$ such that each corresponding value function $u^m$ satisfies Assumption \ref{assum:smooth},
		and $(\mu^m, \sigma^m, g^m) \longrightarrow (\mu, \sigma, g)$ pointwisely,
		one can easily deduce that the corresponding estimators $\psih^m \longrightarrow \psih$ pointwisely.
		If one can check in addition that $(\psih^m)_{m \ge 1}$ is uniformly integrable,
		then  \eqref{eq:unbiased_estim}    still holds,
		see Corollary \ref{coro:approx_repres} below.
	\end{Remark}

\subsection{Proof of Theorem \ref{theo:renorm_general_vol}}
\label{subsec:proof_mainthm}

	{The proof of} Theorem \ref{theo:renorm_general_vol} {is based on} the path-dependent PDE  satisfied by the value {function $u$ associated to }$V_0$ {in \eqref{eq:def_u}}.
	{Namely, let} $\D([0,T])$ denote the Skorokhod space of all $\R^d$-valued c\`adl\`ag paths on $[0,T]$ and {define}
	\begin{equation} \label{eq:def_ub}
		\bar u(t, \om) 
		~:=~
		\E \big[ g \big(X_T, I_T \big) \big| X_{t \wedge \cdot} = \om_{t \wedge \cdot} \big],\; {(t,\om) \in [0,T]\x \D([0,T])}.
	\end{equation}
	Then,
	$$
		\bar u(t, \om) 
		~=~
		u \big( t, \om_t, I_t(\om) \big),
		~\mbox{with}~
		I_t(\om) ~:= \int_0^t \om_s dA_s,
		~~\mbox{for all}~
		(t,\om) \in [0,T]\x \D([0,T]).
	$$
	Let us also recall {the definition of} Dupire's  horizontal and vertical derivatives of the path-dependent functional $\bar u: [0,T] \x \D([0,T]) \longrightarrow \R$, {see \cite{Dupire}}:
	$$
		\partial_t \bar u(t, \om) := \lim_{\eps \searrow 0} \frac{\bar u(t+ \eps, \om_{t \wedge \cdot}) - u(t, \om)}{\eps},
	$$
	and
	\begin{equation} \label{eq:def_uom}
		\partial_{\om} \bar u(t, \om) := \lim_{y \to 0} \frac{\bar u(t, \om + y \1_{\{\cdot \ge t\}}) - u(t, \om)}{y}.
	\end{equation}
	Similarly, one can defined the second order vertical derivative $\partial^2_{\om \om} \bar u$, {whenever these quantities are well-defined, which is the case under Assumption \ref{assum:smooth} since}
	$$
		\partial_{\om} \bar u(t, \om) = D_x u(t, x, \xb)
		~~~\mbox{and}~~
		\partial^2_{\om\om} \bar u(t,\om) = D^2_{xx} u(t, x, \xb),
	$$
	with
	$$
		(x, \xb) := \Big( \om_t, \int_0^t \om_s dA_s \Big).
	$$
	{Then},  {it follows from} Cont and Fourni\'e \cite{ContFournie} that $\bar u$ satisfies the linear path-dependent PDE
	\begin{equation} \label{eq:PPDE_ub}
		\partial_t \bar u(t, \om) + \mu(t, \om_t, I_t(\om)) \cdot \partial_{\om} \bar u(t, \om) + \frac12 \mathrm{Tr} \big[ \sigma \sigma^{\top}(t, \om_t, I_t(\om)) \partial^2_{\om\om} \bar u(t, \om) \big] = 0,\; 
	\end{equation}
	at any ${(t,\om) \in [0,T)\x \D([0,T])}$, with terminal condition $\bar u(T, \om) = g(\om_T, I_T(\om))$, {$\om  \in  \D([0,T])$.}

	\vspace{0.5em}

	For the ease of notation, let us also define, for all $s < t$,
	$$
		\overline A_{s,t} := \frac{1}{t-s} \int_s^t A_{r} d r,
		~~~~
		m^1_{s,t} := \frac{1}{t-s} \int_s^t (A_{r} - A_s) d r,
	$$
	and 
	$$
		m^2_{s,t} := \frac{1}{t-s} \int_0^t \big(A_{r} - \overline A_{s,t} \big)^2 d r,
		~~~~
		\mt^2_{s,t} := \frac{1}{t-s}\int_s^t (A_{r} - A_s)^2 d r.
	$$

	\vspace{0.5em}

	\noindent {\bf Proof of Theorem \ref{theo:renorm_general_vol}.}
	By symmetry of the Brownian motion, $\psit_n$ and $\psit$ have the same distribution as $\psih_n$ and $\psih$. We therefore only consider   $\psih_n$ and $\psih$ in the following.
	
	\vspace{0.5em}
	
	\noindent $\mathrm{(i)}$ Let us first consider the case $n =1$. 
	Recall that $V_0 = \bar u(0, \om)$ where $\bar u$ satisfies the PPDE \eqref{eq:PPDE_ub} and $\om_0 = x_0$.
	Given fixed constants $\mu_0 \in \R^d$ and $\sigma_0 \in \S^d$, we can rewrite the PPDE \eqref{eq:PPDE_ub} as 
	\begin{align*}
		\partial_t&  \bar u(t, \om) 
		+ \mu_0 \cdot \partial_{\om} \bar u(t, \om) + \frac12 \mathrm{Tr} \big[ \sigma_0 \sigma^{\top}_0 \partial^2_{\om\om} \bar u(t, \om) \big]
		+ \bar f(t, \om)
		=0,
	\end{align*}
	with
	$$
		\bar f(t, \om) := f(t, \om_t, I_t(\om)),
	$$
	and
	$$
		f(t, x, \xb) 
		:=
		(\mu(t, x, \xb) - \mu_0)  \cdot D_x u(t, x, \xb) + \frac12 \mathrm{Tr} \big[ \big(\sigma \sigma^{\top} (t, x, \xb) - \sigma_0 \sigma^{\top}_0 \big) D^2_{x} u(t, x, \xb) \big].
	$$
	Then, by the Feynman-Kac's formula in \cite{ContFournie},
	\begin{equation} \label{eq:repres_ub1}
		\bar u(t, \om) 
		~=~
		u (t, \om_t, I_t(\om)) 
		~=~
		\E\Big[
			g \big( \Xb^{t,\om}_T, I_T(\Xb^{t,\om})  \big)
			+ \int_t^T f \big(\Xb^{t,\om}_s, I_s( \Xb^{t,\om}) \big) ds
		\Big],
	\end{equation}
	where 
	$$
		\Xb^{t,\omega}_s = \om_s \1_{\{s < t\}} + \Big( \om_t + \mu_0 (s-t) + \sigma_0 (W_s - W_t) \Big) \1_{\{s \ge t\}}.
	$$
	Notice that $\mu_0$ and $\sigma_0$ can be chosen arbitrary.
	By letting $t=0$, $\om_0 = x_0$, $\mu_0 = \mu(0,x_0, 0)$ and $\sigma_0 = \sigma(0, x_0, 0)$, 
	and recalling that $T_1 = \tau_1 \wedge T$ where $\tau_1$ is a random variable independent of $W$ {with} density function $\rho$ and distribution function $F$,
	it follows that
	\begin{align*}
		V_0 
		&~=~ \E\Big[
			g \big( \Xb^{0,x_0}_T, I_T(\Xb^{0,x_0}) \big)
			+ \int_0^T f \big(s, \Xb^{0, x_0}_s, I_s( \Xb^{0, x_0}) \big) ds
		\Big] \\
		&~=~
		\E \Big[ 
			\frac{g \big(\Xh_{{T}}, I_{{T}}(\Xh) \big)}{1 - F({T})} \1_{\{T_1 \ge T\}} 
			+ 
			\frac{f \big(T_1, \Xh_{T_1}, I_{T_1}(\Xh) \big)}{\rho(T_1)} \1_{\{T_1 < T\}} 
		\Big]
		~=~
		\E \big[ \psih_1 \big].
	\end{align*}

	\noindent $\mathrm{(ii)}$ Next, {define for $s<t$ and $\omega\in \D([0,T])$}
	\begin{equation} \label{eq:def_Mst}
		M_{s,t}^{{\omega}}
		~:=
		\frac{1}{(t-s) m^2_{s,t}} ({\sigma_s^{\omega \top}})^{-1} \int_s^t \Big( \mt^2_{s,t} - m^1_{s,t} (A_u - A_s) \Big) dW_u
	\end{equation}
	with
	$$
	\sigma_s^{\omega}:=\sigma{(s,\om_s,I_s(\om))}.
	$$
	Also define 
	\begin{align*}
	f^{\omega}_{s,t}(x, \xb):=&(\mu(t, x, \xb) - \mu(s,\omega_s,I_s(\omega))  \cdot D_x u(t, x, \xb) \\
	&+ \frac12 \mathrm{Tr} \big[ \big(\sigma \sigma^{\top} (t, x, \xb) - \sigma \sigma^{\top}(s,\omega_s,I_s(\omega) \big) D^2_{x} u(t, x, \xb) \big].
	\end{align*}
	In view of \eqref{eq:repres_ub1} and Lemma \ref{lemm:auto_diff} (see also {\eqref{eq:def_Ih1t}}), one has, {for any $n$},
	\begin{align} \label{eq:repres_Dub1}
		& \partial_{\om} \bar u(t,\om) 
		~=~
		D_x u(t, \om_t, I_t(\om)) \nonumber \\
		=~& 
		\E \Big[
			g \big( \Xb^{t,\om}_T , I_T(\Xb^{t,\om})  \big) M_{t,T}^{\omega} + \int_t^T f{_{t,s}^\omega}\big(\Xb^{t,\om}_s, I_s( \Xb^{t,\om}) \big) M_{t,s}^{{\omega}} ds
		\Big] \nonumber \\
		=~& 
		\E \Big[
			\Big( g \big( \Xb^{t,\om}_T, I_T(\Xb^{t,\om}) \big) - g(\om_t, I_t(\om)) \Big) M_{t,T}^{\omega} + \int_t^T \!\! f{_{t,s}^\omega} \big(\Xb^{t,\om}_s, I_s( \Xb^{t,\om}) \big) M_{t,s}^{{\omega}} ds
		\Big] \nonumber \\
		=~&
		\E \Big[ 
			\frac{g \big( \Xb^{t,\om}_T, I_T(\Xb^{t,\om}) \big) - g(\om_t, I_t(\om))}{1 - F(T-t)} ~M_{t,T}^{{\omega}}~ \1_{\{t+\tau_n \ge T\}} \nonumber \\
		&~~~~~~~~~~~~~~~~~~
			+ 
			\frac{ f{_{t,t+ \tau_n}^\omega} \big(\Xb^{t,\om}_{t+ \tau_n}, I_{t+\tau_n}(\Xb^{t,\om}) \big)}{\rho(\tau_n)} ~M_{t, t+\tau_n}^{{\omega}} ~\1_{\{t+\tau_n < T\}} 
		\Big],
	\end{align}
	{in which we used that $\E[M_{t, T}^{\omega}]=0$ and that   $\tau_n$ is independent of $W$ and has $\rho$ for density. Similarly, for any $n$,}  
	\begin{align}  \label{eq:repres_D2ub1}
		&\partial^2_{\om \om} \bar u(t,\om) ~=~ D^2_{xx} u(t, \om_t, I_t(\om)) \nonumber \\
		=~& 
		\E \Big[
			\Big( g \big( \Xb^{t,\om}_T, I_T(\Xb^{t,\om}) \big) - g(\om_t, I_t(\om)) \Big) \Big( M_{t,T}^{{\omega}} M_{t,T}^{{{\omega}}\top} - \frac{1}{T-t} \frac{\mt^2_{t,T} }{m^2_{t,T}} \big({\sigma_t^{\omega } \sigma_t^{\omega \top})} \big)^{-1} \Big) \nonumber \\
			&~~~~~~~~~~~~~~~~~+ \int_t^T f_{t,s}^\omega \big(\Xb^{t,\om}_s, I_s( \Xb^{t,\om}) \big) \Big( M_{t,s}^{{\omega}} M_{t,s}^{{{\omega}}\top} - \frac{1}{t-s} \frac{\mt^2_{t,s} }{m^2_{s,t}} \big(\sigma_t^{\omega} \sigma_t^{\omega \top}\big)^{-1} \Big) ds
		\Big] \nonumber \\
		=~& 
		\E \bigg[
			\frac{ g \big( \Xb^{t,\om}_T, I_T(\Xb^{t,\om}) \big) - g(\om_t, I_t(\om))}{1-F(T-t)} \Big( M^{{\omega}}_{t,T} M_{t,T}^{{{\omega}}\top} - \frac{1}{T-t} \frac{\mt^2_{t,T} }{m^2_{t,T}} \big(\sigma_t^{\omega} \sigma_t^{\omega \top}\big)^{-1} \Big) 
				\1_{\{t+\tau_n \ge T\}} \nonumber \\
			&~~+~ \frac{ f_{t,t+\tau_n}^\omega \big(\Xb^{t,\om}_{t+\tau_n}, I_{t+\tau_n}( \Xb^{t,\om}) \big)}{\rho(\tau_n)}
				 \Big( M^{{\omega}}_{t,t+\tau_n} M_{t,t+\tau_n}^{{{\omega}}\top} - \frac{1}{\tau_n} \frac{\mt^2_{t,t+\tau_n} }{m^2_{t, t+ \tau_n}} \big(\sigma_t^{\omega} \sigma_t^{\omega \top}\big)^{-1} \Big) \1_{\{t+\tau_n < T\}}
		\bigg].
	\end{align}
	
	\noindent $\mathrm{(iii)}$ Take now $n\ge 2$ and observe that 	 $\psih_{n-1}$ involves  $\partial_{\om} \bar u$ and $\partial^2_{\om\om} \bar u$ (or equivalently $D_x u$ and $D^2_{xx} u$) through the definition of $f$.
	By setting $t = T_{n-1}$, $\om_{T_{n-1} \wedge \cdot} = \Xh_{T_{n-1} \wedge \cdot}$  in \eqref{eq:repres_Dub1} and \eqref{eq:repres_D2ub1},
	and then plugging the representation formula of $\partial_{\om} \bar u (T_{n-1}, \Xh)$ and $\partial^2_{\om\om} \bar u(T_{n-1}, \Xh)$ into $\psih_{n-1}$ in \eqref{eq:def_psi_n},
	it follows that 
	$$
		\E \big[ \psih_{n-1} \big] 
		~=~
		\E \big[ \psih_{n} \big].
	$$
	Since $V_0=\E[\psih_1]$ by (i), it follows that  $V_0=\E[\psih_{n}]$ for all $n\ge 1$ by induction.

	\vspace{0.5em}

	\noindent $\mathrm{(iv)}$ Finally, we notice that
	$$
		\lim_{n \to \infty} \psih_n 
		~=~
		\psih,
		~~\mbox{a.s.}
	$$
	{Since we assumed that  $(\psih_n)_{n \ge 1}$ is uniformly integrable,}
	it follows that $\psih$ is integrable and satisfies 
	\b*
		V_0 ~~=~~ \lim_{n \to \infty} \E \big[  \psih_n \big] 
		~~=~~  \E \big[  \lim_{n \to \infty} \psih_n \big]  ~~=~~ \E \big[ \psih \big].
	\e*
	\qed

\section{Integrability under structural conditions on the coefficients}
\label{sec:integrability}

	Let us provide here some sufficient conditions to ensure the (square) integrability conditions of $\psi$, required in Theorem \ref{theo:renorm_general_vol}.

\subsection{An upper bound of the estimators}

	Recall that $A$ is a (deterministic) function with finite variation and let $|A|$ be the corresponding total variation process. 
	We now impose a structural condition relating the path-regularity of $|A|$ and the regularity of the coefficients of the diffusion process $X$, in the spirit of \cite{BouchardTan}.
	Recall also that a modulus of continuity $\varpi$ is an increasing real-extended valued function $\varpi: [0, \infty] \longrightarrow [0, \infty]$ such that $\lim_{t \searrow 0} \varpi(t) = 0$.
		
	\begin{Assumption} \label{assum:mu_sigma}
		There are constants $L > 0$ and $\alpha_1, ~\alpha_2 \in (0, 1]$,
		together with modulus of continuity $\varpi_1 (\cdot)$ and $\varpi_2(\cdot)$, such that 
		\begin{equation} \label{eq:cond_varpi}
			\varpi_{i} \big( C  \big( |A|_t - |A|_s \big) \big)
			\le
			L ~C~ (t-s)^{\alpha_i/2},
			~\mbox{for all}~
			0 \le s \le t \le T, ~ C \ge 1,
			~i= 1, 2.
		\end{equation}

		\noindent $\mathrm{(i)}$ The functions $(\mu, \sigma)  : [0,T] \x \R^d \x \R^d \to \R^d \x \S^d$ are uniformly bounded,
		and for all $ (t,x,\xb), ~(s,y, \bar y) \in [0,T] \x \R^d \x \R^d$,
		\begin{equation} \label{eq:Lipschitz_coef2}
			\big| \big(\mu, \sigma \big) (t,x, \xb) - \big( \mu, \sigma\big) (s,y, \bar y) \big|
			~\le~
			L \Big( |t-s|^{\alpha_1/2} + \big| x-y \big|^{\alpha_1} + \varpi_1 \big( | \xb - \bar y | \big) \Big).
		\end{equation}
		
		\noindent $\mathrm{(ii)}$ For all $(x, \xb), ~(y,  \bar y) \in \R^d \x \R^d$,
		\begin{equation} \label{eq:cond_g}
			\big| g(x, \bar x) - g(y, \bar y) \big|
			~\le~
			L  \Big( \big| x - y \big|^{\alpha_2} + \varpi_2 \big( | \bar x - \bar y | \big) \Big).
		\end{equation}
		
		\noindent $\mathrm{(iii)}$ For all $(x, \xb), ~(\Delta x, \Delta \bar x) \in \R^d \x \R^d$,
		\begin{equation} \label{eq:cond_antithetique}
			\Big| g\big(x+\Delta x, \bar x + \Delta \bar x  \big) + g\big(x- \Delta x, \bar x + \Delta \bar x \big) - 2 g\big(x, \bar x\big) \Big| 
			~\le~
			L \Big( \big| \Delta x \big|^{2 \alpha_2}  +  \varpi_2 \big( | \Delta \bar x | \big)^2 \Big).
		\end{equation}		
	\end{Assumption}

	\begin{Remark}
		$\mathrm{(i)}$ Assume that there is some constant $\alpha_0 \ge \frac{\alpha_1}{2}$ such that $A$ satisfies $ |A|_t - |A|_s \le (t-s)^{\alpha_0}$ for all $0 \le s \le t \le T$,
		and $(\mu, \sigma)(t, x, \bar x)$ is $\frac{\alpha_1}{2}$-H\"older in $t$, $\alpha_1$-H\"older in $x$ and $\frac{\alpha_1}{2 \alpha_0}$-H\"older in $\bar x$,
		then \eqref{eq:Lipschitz_coef2} holds true with $\varpi_1(s) := s^{\alpha_1/(2 \alpha_0)}$ for all $s \ge 0$.

		\vspace{0.5em}

		Similarly, if $g(x, \bar x)$ is $\alpha_2$-H\"older in $x$ and $\frac{\alpha_2}{2 \alpha_0}$-H\"older in $\bar x$,
		then \eqref{eq:cond_g} holds true with $\varpi_2(s) := s^{\alpha_1/(2 \alpha_0)}$ for all $s \ge 0$.

		\vspace{0.5em}
		
		\noindent $\mathrm{(ii)}$ The condition \eqref{eq:cond_antithetique}  is equivalent to
		$$
			\Big| g\big(x, \bar x + \Delta \bar x  \big) -  g\big(x, \bar x\big) \Big| 
			~\le~
			L \varpi_2 \big( | \Delta \bar x | \big)^2,
		$$
		and
		$$
			\Big| g\big(x+\Delta x, \bar x \big) + g\big(x- \Delta x, \bar x  \big) - 2 g\big(x, \bar x\big) \Big| 
			~\le~
			L \big| \Delta x \big|^{2 \alpha_2} .
		$$

		\vspace{0.5em}

		Let $\alpha_2  \in (\frac12,1]$, $(A_t)_{t \in [0,T]}$ satisfies $ |A|_t - |A|_s \le (t-s)^{\alpha_0}$ for some $\alpha_0 \ge {\alpha_2}$.
		Assume that
		$g(x, \bar x)$ is continuous differentiable in $x$, $D_x g(x, \bar x)$ is $(2\alpha_2-1)$-H\"older in $x$, and $g(x, \bar x)$ is $\frac{\alpha_2}{\alpha_0}$-H\"older in $\bar x$,
		then \eqref{eq:cond_antithetique} holds true with $\varpi_2(s) := s^{\alpha_2/(2 \alpha_0)}$ for all $s \ge 0$.

	\end{Remark}

	Under this condition, each moments of $|\psi_n|$ can be controlled by the norm of the simpler quantity, which does not depend on $(\mu,\sigma)$ anymore, 
	\begin{align*}
		\phi^C_n 
		~:=~&
			\1_{\{N_T \le n-1\}} ~\frac{C \Delta T_{N_T+1}^{\alpha_2} }{1 - F(\Delta T_{N_T+1})} 
				\prod_{k=1}^{N_T} \Big( 
					\frac{C \Delta T_{k}^{\alpha_1/2}  }{\rho(\Delta T_k) }  \frac{1}{\Delta T_{k+1}} \frac{\mt^2_{k,k+1}}{m^2_{k, k+1}} 
				\Big) \\
			&~~~~~~~~~+~ 
			\1_{\{N_T \ge n\}} ~ \frac{C \Delta T_{n}^{\alpha_1/2}}{~\rho(\Delta T_n)~}~ 
				\prod_{k=1}^{n-1}  \Big( 
					\frac{C \Delta T_{k}^{\alpha_1/2}  }{\rho(\Delta T_k) }  \frac{1}{\Delta T_{k+1}} \frac{\mt^2_{k,k+1}}{m^2_{k, k+1}} 
				\Big)
	\end{align*}
	and moments of  $|\psi|$ can be controlled through
	$$
		\phi^C ~:=~ \lim_{n \to \infty} \phi^C_n 
		~=~
		\frac{C \Delta T_{N_T+1}^{\alpha_2}}{1 - F(\Delta T_{N_T+1})} 
				\prod_{k=1}^{N_T} \Big( 
					\frac{C \Delta T_{k}^{\alpha_1/2}  }{\rho(\Delta T_k) }  \frac{1}{\Delta T_{k+1}} \frac{\mt^2_{k,k+1}}{m^2_{k, k+1}} 
				\Big),
	$$
	for a suitable constants $C$. 
	To state this result, we need to introduce the sub-$\sigma$-field
	$$
		\Gc_{\T} := \sigma(\tau_k~: k \ge 1).
	$$
	\begin{Proposition} \label{prop:bound_psih}
		Let Assumptions \ref{assum:smooth} and \ref{assum:mu_sigma} hold true.
		Then, for each  $p \ge 1$, there exist a constant $C_{1,p} > 0$ which depends only on $p \ge 1$, $\eps_0 > 0$ in Assumption \ref{assum:smooth},  $L >0$ in Assumption \ref{assum:mu_sigma} and the $L^{\infty}$-norm of $(D_x u, D^2_{xx} u, g)$, such that
		$$
			\E \Big[ \big| \psi_n \big|^p  ~\Big| \Gc_{\T} \Big]
			~\le~
			\big| \phi^{C_{1,p}}_n \big|^p, ~\mbox{a.s.},
			~~\mbox{for all}~ n \ge 1,
		$$
		and a constant $C_{2,p} > 0$ which  depends only on the constants $p \ge 1$, $\eps_0 > 0$ in Assumption \ref{assum:smooth} and  $L >0$ in Assumption \ref{assum:mu_sigma},  such that
		\begin{equation} \label{eq:integ_psi_p}
			\E \Big[ \big| \psi \big|^p  ~\Big| \Gc_{\T} \Big]
			~\le~
			\big| \phi^{C_{2,p}} \big|^p,
			~\mbox{a.s.}
		\end{equation}
	\end{Proposition}
	\proof
		First, we notice that, for some constant $C > 0$ depending only on the constant $L >0$ in Assumption \ref{assum:mu_sigma} and $\varepsilon_0>0$ in Assumption \ref{assum:smooth},
		$$
			\Big| \Wch^1_k \big|
			~\le~
			C \Big( \Delta T_{k}^{\alpha_1/2} + \big| \Xh_{T_{k}} - \Xh_{T_{k-1}} \big|^{\alpha_1} + \varpi_{1} \big( \big| \Ih_{T_{k}}- \Ih_{T_{k-1}} \big| \big) \Big) 
			\big| M_{k+1} \big| ,
		$$
		and
		$$
			\Big| \Wch^2_k \big|
			~\le~
			C \Big( \Delta T_{k}^{\alpha_1/2} + \big| \Xh_{T_{k}} - \Xh_{T_{k-1}} \big|^{\alpha_1} + \varpi_{1} \big( \big| \Ih_{T_{k}}- \Ih_{T_{k-1}} \big| \big) \Big) 
			\Big( \big| M_{k+1} \big|^2 +  \frac{1}{\Delta T_{k+1}} \frac{\mt^2_{k,k+1}}{m^2_{k,k+1}} \Big).
		$$
		Further, there is some constant $C> 0$ such that
		\begin{align*}
			&
			\big| g(\Xh_T, \widehat I_T) - g(\Xt_{T}, \widetilde I_{T}) \big| \\
			~\le~&
			C \Big( \big| \Delta T_{N_T+1} \big|^{\alpha_2} + \big| \sigma_{T_{N_T}} \Delta W_{N_T+1} \big|^{\alpha_2} + \varpi_{2} \big( \big| \Ih_{T}- \Ih_{T_{N_T}} \big| \big)+ \varpi_{2} \big( \big| \widetilde I_{T}- \Ih_{T_{N_T}} \big| \big) \Big),
		\end{align*}
		and
		\begin{align*}
			&
			\big| g(\Xb_T, \overline I_T) + g(\widetilde X_T, \widetilde I_T) - 2 g(\Xb_{T_{N_T}}, \overline I_{T_{N_T}}) \big| \\
			\le~&
			C\Big( \big| \Delta T_{N_T+1} \big|^{2 \alpha_2} + \big|  \sigma_{T_{N_T}}  \Delta W_{N_T+1} \big|^{2 \alpha_2} 
				+ \varpi_{2} \big( \big| \Ih_{T}- \Ih_{T_{N_T}} \big| \big)^2
				+ \varpi_{2} \big( \big| \widetilde I_{T}- \Ih_{T_{N_T}} \big| \big)^2
			\Big).	
		\end{align*}
		Notice that, for any $q > 0$, there exists $C_q > 0$ such that
		$$
			\E \Big[ \big| \Delta W_k \big|^{2q} +  \big| \Xh_{T_{k}} - \Xh_{T_{k-1}} \big|^{2q} \Big| \Gc_{\T} \Big]
			~\le~
			C_q \Delta T_k^q,
		$$
		and by \eqref{eq:cond_varpi},
		$$
			\E \Big[ \varpi_i \big( \Ih_{T_k} - \Ih_{T_{k-1}} \big)^{q} \Big| \Gc_{\T} \Big] 
			~\le~
			C \E \Big[ \max_{0 \le t \le T} |X_t|^q \Big] \Delta T_k^{\alpha_i q /2}
			~\le~
			C_q \Delta T_k^{\alpha_i q/2},
			~i=1,2.
		$$
		Together with Lemma \ref{lem:Moment_estimate} and the definition of $\psi_n$ and $\psi$ in and above \eqref{eq:def_psi},
		it follows that, for any $p \ge 1$,
		$$
			\E \Big[ \big| \psi_n \big|^p  ~\Big| \Gc_{\T} \Big]
			~\le~
			\big| \phi^{C_{1,p}}_n \big|^p, ~\mbox{a.s.},
		$$
		and
		$$
			\E \Big[ \big| \psi \big|^p  ~\Big| \Gc_{\T} \Big]
			~\le~
			\big| \phi^{C_{2,p}} \big|^p,
			~\mbox{a.s.}
		$$
		for a constant $C_{1,p} > 0$ depending only on $\eps_0$, $L$ and the $L^{\infty}$ norm of $D_x u$, $D^2_{xx} u$ and $g$,
		and a constant $C_{2,p} > 0$ depending only on $\eps_0$ and $L$.
	\qed

	\begin{Remark} \label{rem:antithetique}
		The estimator $\psi$ in \eqref{eq:def_psi}  has a better integrability property than the estimator $\psih$ in \eqref{eq:def_psih} in general.
		Let us consider for example the case where
		$$
			A_t  = t,~~\mbox{for all}~~
			t \in [0,T],
			~~\mbox{so that}~~
			\frac{\mt^2_{k,k+1}}{m^2_{k,k+1}}
			= 4
			~\mbox{(see Lemma \ref{lem:Moment_estimate})}.
		$$

		For $\psih$, one can naturally apply \eqref{eq:cond_g} to obtain the bound
		$$
			\E \Big[ \big| \big( g(\Xh_T, \Ih_T) - g(\Xh_{T_{N_T}}, \Ih_{T_{N_T}}) \big)  \Wch^2_{N_T} \big| ~\Big|~ \Gc_{\T} \Big]
			~\le~
			C \Delta T_{N_T+1}^{- 1 + \alpha_2/2},
		$$ 
		where the r.h.s. term is generally not square integrable for $\alpha_2 \in (0, 1]$ (see also Lemma \ref{lemm:MittagLeffler} below).
		
		\vspace{0.5em}
		
		However, for the estimator $\psi$, one can apply  \eqref{eq:cond_antithetique} to obtain the bound
		$$
			\E \Big[ \big| \big( g(\Xh_T, \Ih_T) + g(\widetilde X_T, \widetilde I_T) - 2 g(\Xh_{T_{N_T}}, \Ih_{T_{N_T}}) \big)  \Wch^2_{N_T} \big| ~\Big|~ \Gc_{\T} \Big]
			~\le~
			C \Delta T_{N_T+1}^{-1 + \alpha_2},
		$$
		where the r.h.s. term becomes square integrable when $\alpha_2 \in ( \frac12, 1]$.
	\end{Remark}

\subsection{A sufficient condition for square integrability}
\label{subsec:sufficient_cond}
 
	We provide here a sufficient condition to ensure the (square) integrability of $\phi^C_n$ and $\phi^C$.
	Recall that $T_k$ is defined in \eqref{eq:def_Tk} with the sequence of i.i.d. random variables $(\tau_k)_{k \ge 1}$, where $\tau_1$ has the density function $\rho$.

	\begin{Theorem} \label{thm:Squ_Integ}
		Let Assumptions \ref{assum:smooth} and \ref{assum:mu_sigma} hold true.
		Suppose that the density function $\rho$ satisfies, for some constants $C_1 > 0$, $C_2 > 0$ and $\kappa_1 > 0$, $\kappa_2 > 0$,
		\begin{equation} \label{eq:cond_rho}
			C_1 s^{\kappa_1 -1}
			~\le~
			\rho(s) 
			~\le~
			C_2 s^{\kappa_2 -1},
			~~\mbox{for all}~
			s \in [0,T].
		\end{equation}
		Suppose further that there exists   $C_3 > 0$ and $\beta \in (0,1]$ such that		
		\begin{equation} \label{eq:m_mt_order}
			m^2_{s,t} > 0,
			~~
			\frac{\mt^2_{s,t}}{m^2_{s,t}} ~\le~ C_3 (t-s)^{- \beta}, ~~\mbox{for all}~ s < t.
		\end{equation}
		Assume in addition that, for some $p \ge 1$,
		\begin{equation} \label{eq:cond_kappa} 
			\kappa_2 + p \Big(\frac{\alpha_1}{2} - \kappa_1 - \beta \Big) > 0
			~~\mbox{and}~~
			p \big(\alpha_2 - \beta - 1 \big) + 1 > 0.
		\end{equation}
		Then for any constant $C> 0$, the sequence $\big( \big| \phi^C_n \big|^p \big)_{n \ge 1}$ is uniformly integrable.
		
		\vspace{0.5em}
		
		Consequently,  if \eqref{eq:cond_kappa} holds for $p=1$, then $\psi$ satisfies the representation \eqref{eq:unbiased_estim}.
		If \eqref{eq:cond_kappa} holds for $p = 2$, then $\psi$ is in addition square integrable.
	\end{Theorem}
	
	\begin{Remark}
	$\mathrm{(i)}$ Let $\rho$ be the density of the Gamma distribution, with parameters $\kappa > 0$ and $\theta > 0$, i.e. 
		$$
			\rho(s) = \frac{s^{\kappa -1} e^{-s/\theta}}{\Gamma(\kappa) \theta^{\kappa}}, ~\mbox{for all}~ s > 0,
		$$
		then it satisfies \eqref{eq:cond_rho} with $\kappa_1 = \kappa_2 = \kappa$.
		
		\vspace{0.5em}
		
	\noindent $\mathrm{(ii)}$
		When $A_t = t$ for all $t \in [0,T]$, it follows by Lemma \ref{lem:Moment_estimate} that 
		$$
			\frac{\mt^2_{s,t}}{m^2_{s,t}} 
			~=~
			4,
			~~\mbox{for all}~ s < t,
		$$
		so that \eqref{eq:m_mt_order} holds true with $\beta = 0$.
		Let us also refer to \cite[Example 2.4]{BouchardTan} for more examples of $A$ satisfying \eqref{eq:m_mt_order} with explicit constant $\beta \ge 0$.
	\end{Remark}

	\begin{Remark}
		Let $A_t = t$ for all $t \in [0,T]$, so that \eqref{eq:m_mt_order} holds true with $\beta = 0$.
		Assume that for some constant $\alpha \in (0,1]$, the functions $(\mu, \sigma)(t,x, \bar x)$ are both $\alpha$-H\"older in $(t, x, \bar x)$, 
		$g(x, \bar x)$ is  $\alpha$-H\"older in $\xb$, and $D_x g(x, \bar x)$ is $\alpha$-H\"older in $x$.
		Then Assumption \ref{assum:mu_sigma} holds with $\alpha_1 = \alpha_2 = \alpha$.
		
		\vspace{0.5em}
	
		\noindent $\mathrm{(i)}$ When $p=1$, the condition \eqref{eq:cond_kappa} holds true as soon as coefficient $\kappa_1, \kappa_2$ in \eqref{eq:cond_rho} satisfies 
		$\kappa_1 - \kappa_2 < \alpha/2$.
	
		\vspace{0.5em}

		\noindent $\mathrm{(ii)}$ When $p=2$,  the condition \eqref{eq:cond_kappa} holds true as soon as $\alpha > \frac12 \vee (2 \kappa_1 - \kappa2)$.
	\end{Remark}

	\noindent {\bf Proof of Theorem \ref{thm:Squ_Integ}.}
		By changing the constant $C > 0$ from line from line, it follows by \eqref{eq:cond_rho} and \eqref{eq:m_mt_order} that
		\begin{align*}
		\phi^C_n 
		~\le~&
			\1_{\{N_T \le n-1\}} ~{ C \Delta T_{N_T+1}^{\alpha_2} }
				\prod_{k=1}^{N_T} \Big( 
					 C \Delta T_{k}^{ \frac{\alpha_1}{2} - \kappa_1 + 1} ~\Delta T_{k+1}^{-(1 + \beta)} 
				\Big) \\
			&~~~~~~~~~+~ 
			\1_{\{N_T \ge n\}} ~ C \Delta T_{n}^{\frac{\alpha_1}{2} - \kappa_1 + 1} ~
				\prod_{k=1}^{n-1}  \Big( 
					 C \Delta T_{k}^{ \frac{\alpha_1}{2} - \kappa_1 + 1} ~\Delta T_{k+1}^{-(1 + \beta)} 
				\Big) \\
		~\le~&
			\1_{\{N_T \le n-1\}} ~ C \Delta T_{N_T+1}^{\alpha_2 - (1+\beta)} 
				\prod_{k=1}^{N_T} \Big( 
					 C \Delta T_{k}^{ \frac{\alpha_1}{2} - \kappa_1 - \beta} 
				\Big) 
			+
			\1_{\{N_T \ge n\}} ~
				\prod_{k=1}^{n}  \Big( 
					 C \Delta T_{k}^{ \frac{\alpha_1}{2} - \kappa_1 - \beta} 
				\Big) .
		\end{align*}
		Therefore, for any $p \ge 1$ and constant $C>0$ big enough, 
		when \eqref{eq:cond_kappa} holds true, or equivalently,
		$$
			\theta := p (\kappa_1 + \beta -  \frac{\alpha_1}{2}) < \kappa_2
			~~\mbox{and}~~
			\eta := 1 + p (\alpha_2 - \beta -1) > 0,
		$$
		it follows by Lemma \ref{lemm:MittagLeffler} that
		$$
			 \E \Big[ \sup_{n \ge 1} \big| \phi^C_n \big|^p \Big] 
			~\le~
			\E \Big[ C^{p+1} \Delta T_{N_T+1}^{\eta - 1}  \prod_{k=1}^{N_T} \Big( C^{p+1} \Delta T_{k}^{ - \theta} \Big)  \Big]
			~<~
			\infty.
		$$
		Consequently, when \eqref{eq:cond_kappa} holds for some $p \ge 1$, 
		the sequence $\big( \big| \phi^C_n \big|^p \big)_{n \ge 1}$ is uniformly integrable for any constant $C>0$.
		
		\vspace{0.5em}
		
		Consequently, it follows by Theorem \ref{theo:renorm_general_vol} and Proposition \ref{prop:bound_psih} that $\psi$ provides the representation result in \eqref{eq:unbiased_estim}.
		When \eqref{eq:cond_kappa} holds for $p = 2$, then by Proposition \ref{prop:bound_psih}, $\psi$ is further square integrable.
	\qed

	\vspace{0.5em}

	Finally, in view of Remark \ref{rem:smoothness}.$\mathrm{(iii)}$ on the smoothness assumption of the value function, we provide here an approximation result.

	\begin{Corollary} \label{coro:approx_repres}
		Let $(\mu^m, \sigma^m, g^m)_{m \ge 1}$ be a sequence of coefficient functions such that $(\mu^m, \sigma^m, g^m) \longrightarrow (\mu, \sigma, g)$ pointwisely.
		Suppose that $(\mu^m, \sigma^m, g^m)_{m \ge 1}$ and $(\mu, \sigma, g)$ all satisfy Assumption \ref{assum:smooth}.$\mathrm{(i)}$ with the same parameters $\eps_0 > 0$ and the same linear growth,
		and all satisfy Assumption \ref{assum:mu_sigma} with uniform parameters $\alpha_1, \alpha_2 \in (0,1]$ and $L > 0$.
		Assume   in addition that, for each $m \ge 1$, the coefficient function $(\mu^m, \sigma^m, g^m)$ induces a value function $u^m$ satisfying the smoothness conditions in Assumption \ref{assum:smooth}.$\mathrm{(ii)}$.
		Then, when condition \eqref{eq:cond_kappa} holds  with $p=2$, one has  
		$$
			\E \big[ \psi^2 \big] < \infty,
			~~\mbox{and}~~
			\E \big[ \psi \big] = \E \big[ g(X_T, I_T) \big] =: V_0.
		$$
	\end{Corollary}
	\proof
		First, for each $m \ge 1$,
		let us denote by $(X^m, I^m)$ the (weak) solution of the SDE \eqref{eq:PSDE} with coefficient functions $(\mu^m, \sigma^m)$, by $\psi^m$ the corresponding representation random variable as defined in \eqref{eq:def_psih}.
		Then it follows by Theorem \ref{thm:Squ_Integ} that 
		$$
			\E \big[ \psi^m \big] = \E \big[ g_m(X^m_T, I^m_T) \big],
			~~\mbox{for each}~~
			m \ge 1.
		$$
		
		Next, by the standard stability results for SDEs, together with the uniqueness of the weak solution $(X, I)$ to the SDE \eqref{eq:PSDE} with coefficient $(\mu, \sigma)$,
		it follows that $(X^m_T, I^m_T) \longrightarrow (X_T, I_T)$ weakly.
		Further, since $(\mu^m, \sigma^m, g^m) \longrightarrow (\mu, \sigma, g)$ pointwisely, one further obtains that
		$$
			\psi^m \longrightarrow \psi,
			~~\mbox{and}~~
			g_m (X^m_T, I^m_T) \longrightarrow g(X_T, I_T),
			~~\mbox{weakly}.
		$$

		Finally, using the growth condition on $g^m$, and the upper bound result in \eqref{eq:integ_psi_p} with uniform constant $C_{2,p}$ for all $m \ge 1$,
		it is easy to deduce that $(\psi^m)_{m \ge 1}$ and $(g_m(X^m_T, I^m_T))_{m \ge 1}$ are both uniformly integrable,
		and that $\psi$ is square integrable.
		It follows then
		$$
			\E \big[ \psi \big] = \E \big[ g(X_T, I_T) \big].
		$$
	\qed

\paragraph{The case with constant volatility coefficient function}
	When $\sigma(\cdot) \equiv \sigma_0$ for some non-degenerate matrix $\sigma_0 \in \S_d$, 
	it is easy to see that $\Wch^2_k = 0$ for all $k \ge 1$.
	In view of Remark \ref{rem:antithetique} and Theorem \ref{thm:Squ_Integ}, one can also consider the variable $\psih$ as a Monte Carlo unbiased estimator.
	In this case, one can drop the conditions in Assumption \ref{assum:mu_sigma}.$\mathrm{(iii)}$ for the (square) integrability analysis of the estimator.

	\vspace{0.5em}

	Since the arguments are almost the same as those in  Proposition \ref{prop:bound_psih} and Theorem \ref{thm:Squ_Integ},
	let us directly provide the representation result without the proof.

	\begin{Proposition}
		Let Assumptions \ref{assum:smooth} and \ref{assum:mu_sigma}.$\mathrm{(i)}$-$\mathrm{(ii)}$ hold true,
		and $\sigma(\cdot) \equiv \sigma_0$ for some non-degenerate matrix $\sigma_0 \in \S_d$.
		Then for each $p \ge 1$, there exists a constant 
		$C_p > 0$
		such that
		\begin{equation} \label{eq:upper_bound_psih}
		\E \Big[ \big| \psih\big|^p \Big| \Gc_{\T}\Big] 
		~\le~
		\left| 
			\frac{C_p \Delta T_{N_T+1}^{\alpha_2/2}}{1 - F(\Delta T_{N_T+1})} 
				\prod_{k=1}^{N_T} \left( 
					\frac{C_p \Delta T_{k}^{\alpha_1/2}  }{\rho(\Delta T_k) } \sqrt{\frac{1}{\Delta T_{k+1}} \frac{\mt^2_{k,k+1}}{m^2_{k, k+1}}} 
				\right)
		\right|^p.
		\end{equation}
		Further, assume that \eqref{eq:cond_rho} holds with parameters $\kappa_1 > 0$, $\kappa_2 > 0$, and \eqref{eq:m_mt_order} holds with parameter $\beta \in (0,1]$,
		and for $p \ge 1$,
		$$
			2 \kappa_2 + p ( 1 + \alpha_1 - 2 \kappa_1 - \beta) > 0 
			~~\mbox{and}~~
			p (\alpha_2 - \beta - 1) + 2 > 0.
		$$
		Then $\big|\psih\big|^p$ is integrable and hence $\psih$ satisfies the representation \eqref{eq:unbiased_estim}.
	\end{Proposition}

	\begin{Remark}
		For the inequality \eqref{eq:upper_bound_psih}, one can easily check that it is enough to choose $C_p > 0$ such that
		$$
			C^p_p 
			~\ge~
			2 L^p \E \Big[ \big ( \big |\Xh_T- \Xh_{T_{N_T}} \big| \big / \Delta T_{N_T+1}^{1/2} \big)^{\alpha_2 p} + L \sup_{0 \le s \le T} |X_s|^p \Big| \Gc_{{T}} \Big]
		$$
		and, for any $k =1, \cdots, N_T$, 
		$$
			C^p_p
			~\ge~
			3 L^p 
			~\E \Big[ \big| \sigma_0^{-1} Z \big|^p \Big]~
			\E \Big[ 1 +   \big( \big |\Xh_{T_k}- \Xh_{T_{k-1}} \big| \big / \Delta T_{k}^{1/2} \big)^{\alpha_1 p}  + L \sup_{0 \le s \le T} |X_s|^p
			 \Big| \Gc_{{T}} \Big],
		$$
		where $Z$ is $d$-dimensional standard Gaussian random vector.
		Thus the constant $C_p > 0$  depends only on the constant $p \ge 1$, the initial condition $x_0$, the Lipschitz constant $L> 0$ of the coefficient functions $(\mu, g)$ in Assumption \ref{assum:mu_sigma} and the norm of $\sigma_0^{-1}$.
		In particular, the variance of the estimator $\psih$ depends directly on the constant $C_p$ with $p=2$.
	\end{Remark}

\section{Numerical examples}
\label{sec:Numerics}

	To illustrate the performance of our estimator, we consider the Asian option pricing problem under two different financial models:
	the Bachelier model and a (path-dependent) local volatility model.

	For our numerical implementation of the unbiased Monte Carlo simulation method, the discrete time grid $(T_k)_{k \ge 1}$ are obtained by simulating the i.i.d. random variables 
	$(\tau_k)_{k \ge 1}$ of density function:
	$$
		\rho(t) ~:=~ \frac{\kappa}{(2T)^{\kappa}} t^{\kappa - 1}  \mathbf{1}_{\{t \in [0,2T]\}},
	$$
	with parameter $\kappa = 0.35$.
	
\subsection{{The} Bachelier model}
	Under the Bachelier model, the underlying risky asset follows the dynamic with constant volatility coefficient:
	\begin{equation} \label{eq:Bachelier}
		dX_t = r X_t dt + \sigma dW_t,
	\end{equation}
	where $r > 0$ is the interest rate, and $W$ is a Brownian motion (under the risk-neutral probability).
	Let us consider the Asian option with payoff
	$$
		(I_T - K)_+, ~~\mbox{where}~I_T := \int_0^T X_t dt,
	$$
	so that the corresponding no-arbitrage price is given by
	$$
		V_0 ~:=~ \E \Big[ e^{-rT} \big(I_T - K \big)_+ \Big].
	$$
	In this setting, one can compute directly that
	$$
		I_T 
		~=
		 \int_0^T X_t dt 
		~=~
		\frac1r \Big( \big(e^{rT} - 1 \big) X_0 + \sigma \int_0^T \big( e^{r(T-t)} - 1 \big) dW_s \Big),
	$$
	and then deduce the reference value:
	$$
		V_0 ~=~ \Sigma ~\Big(  \frac{1}{\sqrt{2 \pi}} e^{- (K- \mu)^2/ (2 \Sigma^2)} -  \frac{K-\mu}{\Sigma} \Phi \big( -  \frac{K-\mu}{\Sigma} \big) \Big),
	$$
	with
	$$
		\mu := \frac1r \big( e^{rT} - 1 \big) X_0,
		~~~
		\Sigma^2 := \frac{\sigma^2}{r^2} \Big( \frac{1}{2r} \big( e^{2rT} - 1 \big) - \frac{2}{r} \big( e^{rT} - 1\big) + T \Big).
	$$

	\vspace{0.5em}
	
	In Table \ref{tab:Bachelier}, we present the numerical simulation results under the Bachelier model.
	For different values of the volatility coefficient $\sigma$, the unbiased estimations are quite close to the reference values.
	Moreover, comparing to the local volatility model (see Tables \ref{tab:LV1} and \ref{tab:LV2}),
	the unbiased simulation method has a better performance under the Bachelier model in terms of  standard deviation.
	Indeed, the term $\Wch^2_k$ in \eqref{eq:Wcb2} has a higher order {of} variance than  $\Wch^1_k$ in \eqref{eq:Wcb1} in general,
	and the term $\Wch^2_k$ disappears in the estimator $\psi$ under the Bachelier model,
	since the volatility coefficient is a constant (see the discussions at the end of Section \ref{subsec:sufficient_cond}).

	\begin{table}[h!]
	\begin{center}
	\begin{tabular}{|l|l|l|l|}
	\hline
		~ & Reference value &  Unbiased MC mean & Std/$\sqrt{N}$     \\ \hline
		$\sigma = 0.05$  &  2.7182  & 2.7159  &  0.0022    \\ \hline
		$\sigma = 0.1$  &  3.6470  & 3.6439  &   0.0034   \\ \hline
		$\sigma = 0.15$  &  4.6960  & 4.6921  &   0.0046   \\ \hline
		$\sigma = 0.2$  &  5.7781  & 5.7733  &   0.0058   \\ \hline
	\end{tabular}
	\end{center}
	\caption{Bachelier model: $r=0.05$, $X_0 = 100$, $K = 100$, $T = 1$, number of simulation $N = 10^7$.}
	\label{tab:Bachelier}
	\end{table}

\subsection{A local volatility model}
	The second model is a path-dependent local volatility model, where the risky asset follows the dynamic
	\begin{equation} \label{eq:BS}
		dX_t = r X_t dt + \sigma(X_t, I_t) dW_t, 
	\end{equation}
	with $\sigma$ taking the values $\sigma_1$ or $\sigma_2$ defined by
	$$
		\sigma_1(x,{\bar x}) = \sigma_0 \big(1+  \sin(x-{\bar x})/4 \big),
		~~\mbox{and}~~
		\sigma_2(x,{\bar x}) = \sigma_0 \big(1 +  \sin(x-{\bar x})/2 \big).
	$$
	We choose the coefficient $r=0.05$, $\sigma_0 = 0.2$, $X_0 = 100$, $K = 100$, $T = 1$.
	For the unbiased simulation method, the number of simulations is $N_1 = 10^8$.
	We also provide the results with the corresponding Euler scheme with different numbers of time steps  $N_T$, and with a number of simulations  $N_2 = 10^6$.
	The results are presented in Tables \ref{tab:LV1} and \ref{tab:LV2}.

	\begin{table}[h!]
	\begin{center}
	\begin{tabular}{|l|c|c|}
	\hline
		 \mbox{Time discretization}  & Euler Scheme mean value &    Std/$\sqrt{N_2}$   \\ \hline
		$N_T = 10$  &  7.04872  & 0.00920483     \\ \hline
		$N_T = 20$  &  7.04872  & 0.00920483     \\ \hline
		$N_T = 40$  &  6.64318  & 0.00871443     \\ \hline
		$N_T = 80$  &  6.57892  & 0.00863623     \\ \hline
		$N_T = 160$  &  6.55219  & 0.00858793     \\ \hline
		$N_T = 200$  &  6.54021  & 0.00856895     \\ \hline \hline
		 ~ &  Unbiased MC mean value &   Std/$\sqrt{N_1}$      \\ \hline
		    & 6.51562  &   0.00951147   \\ \hline 
	\end{tabular}
	\end{center}
	\caption{The local volatility model: $\sigma (x,{\bar x}) = \sigma_1(x,{\bar x})$.}
	\label{tab:LV1}
	\end{table}

	\begin{table}[h!]
	\begin{center}
	\begin{tabular}{|l|c|c|}
	\hline
		 \mbox{Time discretization}  & Euler Scheme mean value &    Std/$\sqrt{N_2}$   \\ \hline
		$N_T = 10$  &  7.87182  & 0.0105428     \\ \hline
		$N_T = 20$  &  7.55318  & 0.0101459     \\ \hline
		$N_T = 40$  &  7.40358  & 0.00997426     \\ \hline
		$N_T = 80$  &  7.32961  & 0.00988389     \\ \hline
		$N_T = 160$  &  7.29838  & 0.00982819     \\ \hline
		$N_T = 200$  &  7.28373  & 0.00980535     \\ \hline \hline
		 ~ &  Unbiased MC mean value &   Std/$\sqrt{N_1}$      \\ \hline
		    & 7.2363  &   0.0484702   \\ \hline 
	\end{tabular}
	\end{center}
	\caption{The local volatility model: $\sigma (x,{\bar x}) = \sigma_2(x,{\bar x})$.}
	\label{tab:LV2}
	\end{table}

	We observe that the results of Euler scheme converges to the unbiased simulation results as the number of time steps $N_T$ increases.
	Moreover, the performance of the unbiased simulation with $\sigma_1$ in Table \ref{tab:LV1} is better than the one with $\sigma_2$ in Table \ref{tab:LV2}.
	The main reason is that $\sigma_1$ has a smaller Lipschitz coefficient than  $\sigma_2$,
	{while} the constant $C_{1,p}$ depends on (and is increasing in) the Lipschitz coefficient $L> 0$ in Proposition \ref{prop:bound_psih}.
	Therefore, {it is not surprising that} the standard deviation of the estimation {associated to} $\sigma_1 $ is smaller that {the one associated to}  $\sigma_2$.

\appendix

\section{Technical lemmas}
\label{sec:technical_lemmas}

	Recall that, for $s < t$, we define
	$$
		\overline A_{s,t} := \frac{1}{t-s} \int_s^t A_{r} d r,
		~~~~
		m^1_{s,t} := \frac{1}{t-s} \int_s^t (A_{r} - A_s) d r,
	$$
	and
	$$
		m^2_{s,t} := \frac{1}{t-s} \int_0^t \big(A_{r} - \overline A_{s,t} \big)^2 d r,
		~~~~
		\mt^2_{s,t} := \frac{1}{t-s}\int_s^t (A_{r} - A_s)^2 d r,
	$$
	and then
	$$
		\Sigma^A_{s,t} 
		~:=~
		\begin{bmatrix}
			(t-s) & \int_s^t (A_t - A_r) dr \\
			\int_s^t (A_t - A_r) dr & \int_s^t (A_t - A_r)^2 dr
		\end{bmatrix},
	$$
	so that
	\begin{align*}
		\mathrm{Det} \big( \Sigma^A_{s,t} \big) 
		&~=~
		(t-s) \int_s^t (A_t - A_r)^2 dr
		- 
		\Big(\!
			\int_s^t (A_t - A_r) dr
		\Big)^2 \\
		&~=~
		(t-s) \int_s^t (A_r - A_s)^2 dr 
		-
		\Big(\!
			\int_s^t (A_r - A_s) dr
		\Big)^2
		= (t-s)^2 m^2_{s,t}.
	\end{align*}
	Further, let us fix a matrix $\sigma_0 \in \S^d$ and a $d$-dimensional standard Brownian motion $W$, we define, for all $x, \bar x \in \R^d$, 
	$$
		\Phi(s, t, x, \bar x, W_{\cdot}) 
		~:=~
		\begin{bmatrix}
			I_d& 0 \\
			(A_t - A_s) I_d & I_d
		\end{bmatrix}
		\begin{bmatrix}
			x \\
			\xb
		\end{bmatrix}
		~+~
		\sigma_0 
		\int_s^t
		\begin{bmatrix}
			I_d \\
			(A_t - A_u) I_d 
		\end{bmatrix}
		dW_u,
	$$
	and (recall  \eqref{eq:def_Mst})
	\begin{equation} \label{eq:def_Mst0}
		M_{s,t} 
		~:=
		\frac{1}{(t-s) m^2_{s,t}} (\sigma^{\top}_0)^{-1} \int_s^t \Big( \mt^2_{s,t} - m^1_{s,t} (A_u - A_s) \Big) dW_u.
	\end{equation}

	\begin{Lemma} \label{lemm:auto_diff}
		Let $g: \R^d \x \R^d \longrightarrow \R$ be a function with polynomial growth, then
		\begin{align*}
			D_x \E \Big[ g\big( \Phi(s, t, x, \xb, W_{\cdot}) \big) \Big] 
			&=
			\E \Big[
				g\big( \Phi(s, t, x, \xb, W_{\cdot}) \big)
				M_{s,t}
			\Big],
		\end{align*}
		and
		\begin{align*}
			D^2_{xx} \E \Big[ g\big( \Phi(s, t, x, \xb, W_{\cdot}) \big) \Big] 
			&=
			\E \Big[
				g\big( \Phi(s, t, x, \xb, W_{\cdot}) \big)
				\Big( 
					M_{s,t} M_{s,t}^{\top} 
					-
					\frac{1}{t-s} \frac{\mt^2_{s,t} }{m^2_{s,t}} \big(\sigma_0 \sigma_0^{\top} \big)^{-1}
				\Big)
			\Big].
		\end{align*}
	\end{Lemma}
	\proof Let us denote the $2 d$-dimensional Gaussian vector
	$$
		Z 
		~:=~
		\sigma_0 
		\int_s^t
		\begin{bmatrix}
			I_d \\
			(A_t - A_u) I_d 
		\end{bmatrix}
		dW_u.
	$$
	Then
	$$
		D_x \E \Big[ g\big( \Phi(s, t, x, \xb, W_{\cdot}) \big) \Big] 
		=
		\begin{bmatrix}
			I_d \\
			 (A_t-A_s) I_d
		\end{bmatrix}^{\top}
		\E \Big[
			g\big( \Phi(s, t, x, \xb, W_{\cdot}) \big)
			\mathrm{Var} \big[ Z \big]^{-1} Z
		\Big],
	$$
	and
	\begin{align*}
		D^2_{xx} \E \Big[ g & \big( \Phi(s, t, x, \xb, W_{\cdot}) \big) \Big]
		=
		\begin{bmatrix}
			I_d \\
			 (A_t-A_s) I_d
		\end{bmatrix}^{\top} \\
		&\E \Big[
			g\big( \Phi(s, t, x, \xb, W_{\cdot}) \big)
			\Big( 
				\mathrm{Var}[Z]^{-1} Z~ \big( \mathrm{Var}[Z]^{-1} Z \big)^{\top} - \mathrm{Var} [Z]^{-1}
			\Big)
		\Big]
		\begin{bmatrix}
			I_d \\
			 (A_t-A_s) I_d
		\end{bmatrix}.
	\end{align*}
	
	By direct computation, one has $\mathrm{Var}[Z] = \sigma_0 \sigma_0^{\top} \Sigma^A_{s,t}$ , 
	so that
	\begin{align*}
		\begin{bmatrix}
			I_d \\
			 (A_t-A_s) I_d
		\end{bmatrix}^{\top}
		\mathrm{Var} \big[ Z \big]^{-1} Z
		=
		\big( \sigma_0 \sigma_0^{\top} \big)^{-1}
		\sigma_0
		\int_s^t 
		\begin{bmatrix}
			1 \\
			A_t-A_s
		\end{bmatrix}^{\top}
		\big( \Sigma^A_{s,t} \big)^{-1} 
		\begin{bmatrix}
			1 \\
			A_t-A_u
		\end{bmatrix}
		d W_u.
	\end{align*}
	Notice that 
	$$
		\big( \Sigma^A_{s,t} \big)^{-1}
		~=~
		\frac{1}{\mathrm{Det} \big( \Sigma^A_{s,t} \big)}
		\begin{bmatrix}
			 \int_s^t (A_t - A_u)^2 du & - \int_s^t (A_t - A_u) du \\
			- \int_s^t (A_t - A_u) du & (t-s)
		\end{bmatrix},
	$$
	it follows by direct computation that
	$$
		\begin{bmatrix}
			I_d \\
			 (A_t-A_s) I_d
		\end{bmatrix}^{\top}
		\mathrm{Var} \big[ Z \big]^{-1} Z
		~=~
		M_{s,t},
	$$
	with $M_{s,t}$ being defined in \eqref{eq:def_Mst0}.	
	
	\vspace{0.5em}
	
	Similarly, one has
	\begin{align*}
		&
		\begin{bmatrix}
			I_d \\
			 (A_t-A_s) I_d
		\end{bmatrix}^{\top} 
		\Big( 
				\mathrm{Var}[Z]^{-1} Z~ \big( \mathrm{Var}[Z]^{-1} Z \big)^{\top} - \mathrm{Var} [Z]^{-1}
			\Big)
		\Big]
		\begin{bmatrix}
			I_d \\
			 (A_t-A_s) I_d
		\end{bmatrix} \\
		=&~
		M_{s,t} M_{s,t}^{\top} 
		-
		\frac{1}{t-s} \frac{\mt^2_{s,t} }{ m^2_{s,t}} \big(\sigma_0 \sigma_0^{\top} \big)^{-1},
	\end{align*}
	which concludes the proof.
	\qed

	\begin{Lemma} \label{lem:Moment_estimate}
		Let $s<t$ and $M_{s,t}$ be defined in \eqref{eq:def_Mst0}.
		Then $M_{s,t}$ is a Gaussian vector, i.e.
		$$
			M_{s,t} 
			~\sim~
			N\Big(0, ~ \frac{1}{t-s} \frac{\mt^2_{s,t} }{m^2_{s,t}} \big(\sigma_0 \sigma_0^{\top} \big)^{-1} \Big).
		$$
		In particular, when $A_t = t$ for all $t \in [0,T]$,
		$$
			\mt^2_{s,t} = \frac{(t-s)^2}{3},
			~~~
			m^2_{s,t} = \frac{(t-s)^2}{12}
			~~\mbox{so that}~~
			\frac{1}{t-s} \frac{\mt^2_{s,t} }{m^2_{s,t}} \big(\sigma_0 \sigma_0^{\top} \big)^{-1}
			=
			\frac{4 }{t-s}  \big( \sigma_0 \sigma_0^{\top} \big)^{-1}.
		$$
	\end{Lemma}
	\proof
	First, it is clear that $M_{s,t}$ is a Gaussian random vector with expectation $0$, and
	\begin{align*}
		\mathrm{Var} \big[ M_{s,t} \big] 
		&~=~
		\frac{1}{ (t-s)^2 \big(  m^2_{s,t} \big)^2} ~  \big(\sigma_0 \sigma_0^{\top} \big)^{-1} \int_s^t \big( \mt^2_{s,t} - m^1_{s,t} ( A_u - A_s) \big)^2 du \\
		&~=~
		\frac{(t-s) \mt^2_{s,t} m^2_{s,t} }{(t-s)^2 \big(  m^2_{s,t} \big)^2 }~  \big(\sigma_0 \sigma_0^{\top} \big)^{-1} 
		~=~
		\frac{1}{t-s} \frac{\mt^2_{s,t}}{m^2_{s,t}} ~  \big(\sigma_0 \sigma_0^{\top} \big)^{-1} 
		.
	\end{align*}
	Next, when $A_t = t$ for all $t \in [0,T]$, it follows by direct computation that
	$$
			\mt^2_{s,t} = \frac{(t-s)^2}{3},
			~~~
			m^2_{s,t} = \frac{(t-s)^2}{12}
			~~\mbox{so that}~~
			\frac{1}{t-s} \frac{\mt^2_{s,t} }{m^2_{s,t}} \big(\sigma_0 \sigma_0^{\top} \big)^{-1}
			=
			\frac{4 }{t-s}  \big( \sigma_0 \sigma_0^{\top} \big)^{-1}.
	$$
	\qed

	Let $T_k$ is defined by \eqref{eq:def_Tk} with the i.i.d. sequence of random variables $(\tau_k)_{k \ge 1}$, where $\tau_1$ has the density function $\rho(s)$.
	\begin{Lemma} \label{lemm:MittagLeffler}
		Assume that the density function $\rho(s)$ satisfies, for some constant $\widetilde C> 0$,
		$$
			\rho(s) \le \widetilde C s^{\kappa - 1},
			~\mbox{for all}~ s \in [0,T].
		$$
		Then for all constants $ C > 0$, $\eta > 0 $ and $\theta < \kappa$, it holds that
		$$
			\E \Big[ \Delta T_{N_T+1}^{\eta-1} \prod_{k=1}^{N_T} \Big( C \Delta T_k^{-\theta} \Big) \Big]
			~<~
			\infty.
		$$
	\end{Lemma}
	\proof
	First, for each $n \ge 0$, we denote $S_n := \{0 < s_1 < \cdots < s_n \}$ and $ds := ds_1 \cdots d s_n$,
	and then compute that
	\begin{align*}
		& \E \Big[ \big(T-  T_{N_T} \big)^{\eta - 1} \prod_{k=1}^{N_T} \Big( C \Delta T_k^{-\theta} \Big) ~\1_{\{N_T = n\}} \Big] \\
		=~&
		\int_{S_n} d s \int_{T-s_n}^{\infty} \rho(t) d t ~(T - s_n)^{\eta - 1} \prod_{k=1}^n \Big( C (s_k - s_{k-1})^{-\theta} \rho(s_k - s_{k-1}) \Big) \\
		\le~ &
		\big( C \widetilde C \big)^n \int_{S_n} \big(T - s_n) \big)^{\eta - 1} \Big( \prod_{k=1}^n (s_k - s_{k-1})^{\kappa - \theta - 1} \Big)  ds \\
		=~&
		\big( C \widetilde C \big)^n  T^{\eta + n( \kappa-\theta) - 1} \frac{\Gamma(\eta) ~\Gamma^n(\kappa - \theta)}{~\Gamma(\eta + n(\kappa - \theta))~}.
	\end{align*}
	Then it follows that, for some constant $K > 0$,
	\begin{align*}
		& \E \Big[ (T- T_{N_T})^{\eta-1} \prod_{k=1}^{N_T} \Big( C \Delta T_k^{-\theta} \Big) \Big]
		\le~
		\sum_{n=0}^{\infty} \frac{K^n}{\Gamma(\eta + n(\kappa - \theta))} 
		~<~ \infty.
	\end{align*}
	\qed

\end{document}